\documentclass{article}
\usepackage{amsmath}
\usepackage{amsfonts,amssymb}
\usepackage[dvips]{graphics}
\usepackage{graphicx}

\input amssym.def

\numberwithin{equation}{section}

\newcommand{\End}{{\rm E n d }}

\newcommand{\Ai}{{A_\infty }}

\newcommand{\Ca}{{\rm {Calc}}}

\newcommand{\Y}{{\bf Y}}

\newcommand{\Ccc}{C_\bullet(A)}

\newcommand{\VM}{\cT^{\bul+1}_{poly}(M)}

\newcommand{\g}{{\mathfrak{g}}^\bullet}

\newcommand{\bul}{{\bullet}}
\newcommand{\gA}{{\mathfrak g}^\bullet_A}
\newcommand{\cT}{{\cal T}}

\newcommand{\CCnA}{{\rm{CC}}^-_{-\bullet}(A)}

\newcommand{\isomoto}{\stackrel{\sim}{\longrightarrow}}
\date{}

\newtheorem{definition}{Definition}

\newtheorem{lemma}{Lemma}
\newtheorem{thm}{Theorem}
\newtheorem{teo}{Theorem}

\newtheorem{pred}{Proposition}
\newtheorem{proposition}{Proposition}

\newtheorem{example}{Example}

\newtheorem{rmk}{Remark}

\title{Noncommutative calculus and the Gauss-Manin connection}

\author{V.A. Dolgushev, D.E. Tamarkin, and B.L. Tsygan}
\date{}

\begin{document}

\large

\maketitle

\begin{center}
{\it To Murray Gerstenhaber on his 80th and to Jim Stasheff on his 70th birthdays}
\end{center}

\begin{abstract}
 After an overview of noncommutative differential calculus, we construct parts of it explicitly and explain why this construction agrees with a fuller version obtained from the theory of operads.
\end{abstract}

%\tableofcontents
\smallskip

{\centerline{\em {}}}

\bigskip

\section{Introduction}\label{s:intro} In this paper we apply the techniques of brace algebras of Gerstenhaber and of $A_\infty$ and $L_{\infty}$ algebras of Stasheff to develop a part of what we call noncommutative differential calculus. Noncommutative calculus is a theory that reconstructs basic algebraic structures arising from the calculus on a manifold in terms of the algebra of functions on this manifold, in a way that works for any algebra, commutative or not. This program is being developed in \cite{DTT1}, \cite{DTT2}, \cite{DTT3}, \cite{TT}, \cite{T}. Let us start by observing that there are several algebraic structures arising from the standard calculus on a manifold:

I. {\bf (Differential graded) Lie algebras and modules over them.} Several key formulas from differential calculus on manifolds, namely the Cartan formulas, use nothing but commutators (the commutator is always understood in the graded sense) and therefore give rise to graded Lie algebras. It has been emphasized by Gelfand and Dorfman \cite{DG} that these
graded Lie algebras and their representations are worthy of being studied and generalized. We use the following notation to describe these algebras.

a) Multivector fields with the Schouten-Nijenhuis-Richardson bracket form a graded Lie algebra that we denote by $\mathfrak{g}^{\bullet}$; $\mathfrak{g}^{\bullet}=\wedge^{\bullet +1}T.$ This graded Lie algebra acts on the space $\Omega ^{-\bullet}$ of forms with reversed grading by the generalized Lie derivative: $L_D=[d,\iota_D]$ where $\iota _D$ is the contraction of a form by a multivector.

b) There is a bigger differential graded Lie algebra (DGLA) $\mathfrak{g}[\epsilon,u]$ with the differential $u\frac{\partial}{\partial \epsilon}$. Here $u$ is a formal parameter of degree $2$ and $\epsilon$ is a formal parameter of degree $1$ such that $\epsilon ^2=0.$ It acts on the complex $\Omega^{-\bullet}[[u]]$ with the differential $ud$ as follows: $X+\epsilon Y$ acts by $L_X+\iota _Y.$

II. {\bf (Differential graded) associative algebras}. There are several:

a) Forms with wedge multiplication.

b) Multivectors with wedge multiplication.

c) Differential operators on functions (and sections of other vector bundles).

d) In particular, differential operators on differential forms.

The algebras IIa) and IIb) are graded commutative.

III. {\bf {Calculi}}. The structures form I, as well as from IIb), give rise to an algebraic structure that we call {\em a calculus}. In particular, multivectors form a {\em Gerstenhaber algebra}, or simply a {\em G-algebra}. We recall the definitions in \ref{ss:c}. They formalize algebraic properties of multivectors with the wedge product and the Schouten bracket, forms with the De Rham differential, and of the former acting on the latter by contraction and by Lie derivative. One reconstructs the algebra IId) from the calculus as an enveloping algebra of a certain kind.

We would like to generalize all the above constructions to the case when a manifold is replaced by a possibly noncommutative algebra over a field of characteristic zero. First note that all the algebras as in I-III will be replaced by {\em algebras up to homotopy}, namely $L_{\infty}$, $A_{\infty}$, $C_\infty$, $G_{\infty}$, and ${\rm{Calc}}_{\infty}$ algebras.

There are two ways to look at such objects. One is to say that a complex $C^{\bullet}$ is an algebra of certain type up to homotopy if a DG algebra $\mathcal{
C}^{\bullet}$ of this type is given, together with a quasi-isomorphism of complexes $\mathcal{
C}^{\bullet}\to C^{\bullet}$. (In case of ${\rm{Calc}}_{\infty}$ algebras one talks rather about pairs of complexes). A $A_\infty$, etc. morphism $C_1\to C_2$ is a chains of morphisms ${\mathcal{C}}_1\leftarrow {\mathcal{C}}\rightarrow {\mathcal{C}}_2$ where the arrow on the left is a quasi-isomorphism. Such a morphism is a {\em quasi-isomorphism} if the map on the right is a quasi-isomorphism, too. There is a natural way of composing such morphisms.

Another way to talk about an algebra up to homotopy is to talk about the complex $C^{\bullet}$ equipped with a series of higher operations satisfying certain relations. We recall the definition of $A_{\infty}$ algebras and modules in terms of higher operations in \ref{ss:aialia}, and (implicitly) an analogous definition of $L_{\infty}$ modules, in the beginning of section \ref{s:tlimsotncc}. The definition of $G_\infty$ algebras in these terms was given in \cite{GJ1}, cf. also \cite{Ta1}, \cite{Hin}; an analogous definition of ${\rm{Calc}}_{\infty}$ algebras will be given in a subsequent work. Morphisms are defined in terms of higher operations as well. In all of the cases discussed above, including the $BV$ case, a homotopy algebra structure can be defined as a coderivation of degree one of a free coalgebra of appropriate type (namely, over the cooperad ${\rm{Lie}}^{\rm{dual}}$, ${\rm{Calc}}^{\rm{dual}}, $ ${\rm{BV}}^{\rm{dual}}$, etc.); that coderivation satisfies the Maurer-Cartan equation. A morphism of two homotopy algebras can be defined as a morphism of resulting DG coalgebras.

One can move between the two ways of defining algebras up to homotopy: from the first definition to the second by a procedure called {\em transfer of structure}, from the second to the first by another procedure called {\em rectification}, cf. \cite{Ka}, \cite{Me}, \cite{KS}, \cite{MV}.

Let $C^{\bullet}(A)$ be the Hochschild cochain complex and $C_{\bullet}(A)$ the Hochschild chain complex of an algebra $A$ over a field of characteristic zero. The former will play the role of noncommutative multivectors and the latter of noncommutative forms. We will start with the noncommutative analog of III and work our way back to I.

It was proven in \cite{DTT2}, \cite{KS} that the pair $C^{\bullet}(A), \;C_{\bullet}(A)$ is a ${\rm{Calc}}_{\infty}$ algebra whose underlying $L_{\infty}$ structure is the one from Hochschild theory, given by the Gerstenhaber bracket on cochains and by some explicit action of cochains on chains. (Those two Lie operations should be viewed as a noncommutative analog of Ia)).

 {\bf Noncommutative version of III.} The ${\rm{Calc}}_{\infty}$ structure from \cite{DTT2} has two related drawbacks: its construction is highly inexplicit and involved, and it is not canonical. The latter part is due to a fundamental fact about ${\rm{Calc}}_{\infty}$ (as well as $G_\infty$) structures: they know how to generate new such structures from themselves. That is, starting from a Gerstenhaber algebra, one can construct new algebras (deformations of the old one) in a universal way, using only the Gerstenhaber operations of the bracket and the product. Similarly, there are universal moves that produce one $G_\infty$ structure from another (indeed, pass from a $G_\infty$ algebra to a DG G-algebra by rectification, then apply the above construction, and go back by the transfer of structure).

 The group generated by these moves is a group of symmetries acting on the set of $G_\infty$ algebra structures on any given space. (Or, which is the same, the operad $G_\infty$ has a nontrivial group of symmetries). It follows from results of \cite{K1}, \cite{KS}, \cite{Tam} that the Grothendieck-Teichm\"{u}ller group maps to this group  of symmetries. Apparently, all of the above is true for ${\rm{Calc}}_{\infty}$ algebras.

This is a phenomenon that is largely absent from the world of associative or Lie algebras (note, however, that a Lie algebra structure $[,\,]$ automatically comes in a one-dimensional family $t[,\,]$). Indeed, higher operations in an $L_\infty$ or an $A_\infty$ algebra are of negative degree, and there is no way to produce a universal formula for such an operation using only the commutator or the product that are of degree zero).

On the other hand, if one has a Gerstenhaber algebra, one can define, say, a new $A_\infty$ structure on it in a universal way, using only the multiplication and the bracket. An example of such a structure over the ring $\mathbb{C}[\hbar]/(\hbar ^2)$: $m_3=0; $ $m_4(a_1,a_2,a_3,a_4)=(-1)^{|a_1|+|a_2|}\hbar[a_1,a_2][a_3,a_4];$  $m_k=0,\;k>4.$ (This $A_\infty$ structure can be extended to a $G_\infty$ structure; to see this, recall from \cite{GJ1}, \cite{Ta1} that a $G_\infty$ structure on $C$ is a collection of operations $m_{k_1,\ldots, k_n}: C^{k_1+\ldots+k_n}\to C$ subject to some relations; let $m_4$ be as above, $m_{2,2}(a_1,a_2;a_3,a_4)=\hbar[[a_1,a_2],[a_3,a_4]]$ and let all other higher operations be zero).

 Note that all of the above applies to the classical calculus structure on the spaces of multivectors and forms; many other $\rm{Calc}_\infty$ structures may be generated from it. But, on the one hand, the standard structure is clearly the most natural; on the other hand, one can prove that all $\rm{Calc}_\infty$ structures that could be written naturally on a smooth manifold are equivalent to the standard one. This can be done using the argument as in \cite{Ta1} and \cite{Hin}, plus some formal differential geometry. (There are nonstandard $\rm{Calc}_\infty$ structures on multivectors and fields; any 3-cohomology class gives one. All this is of course consistent with the fact that there are no natural odd cohomology classes on manifolds; there are natural even classes, namely the Chern classes of the tangent bundle).

 If we try to look at the noncommutative analogs of the structures I and II that arise from the noncommutative version of III, the situation becomes easier for II and a lot easier for I.

 {\bf Noncommutative version of II.} The structures a) does not generalize to our version of noncommutative calculus. Indeed, there is a product on the Hochschild homology, but only for a commutative ring $A$. For example, the degree zero Hochschild homology of $A$ is $A/[A,A]$, the quotient of $A$ by the linear span of commutators; this space does not have any natural multiplication. In comparison, the zero degree Hochschild cohomology of $A$ is the center of $A$ which is always a commutative algebra. (Note, however, that for a deformation quantization of a smooth manifold the space of noncommutative forms, i.e. the Hochschild chain complex, is quasi-isomorphic to the Poisson chain complex; this follows from \cite{D} and \cite{Sh}, cf. also \cite{TT1}. But the differential in the Poisson chain complex is a BV operator. Therefore, in the case of deformation quantization, the Hochschild chain complex is {\em a homotopy BV algebra}; that is, there is a natural model for chains that carries a graded commutative product; the differential is not a derivation with respect to this product but rather a BV operator (in particular, a differential operator of order two). As for the structure c), its generalization to our version of noncommutative calculus is unknown and, in our view, not likely to exist.

 The noncommutative version of the DG algebra IId) of differential operators on differential forms was described in \cite{TT} (we recall and use it in this paper). It was proven there that, indeed, it is the one coming from the $\rm{Calc}_\infty $ structure generalizing III.

 As for the generalized algebra IIb) of multivectors, the situation is more delicate. It is easy to name one candidate, the DG algebra $C^{\bullet}(A)$ of Hochschild cochains with the standard differential $\delta$ and the cup product $\smile$. All we know is that the $C_{\infty}$ algebra structure on $C^{\bullet}(A)$ which is a part of the $\rm{Calc}_\infty $ structure from \cite{DTT2} and \cite{KS} is an $A_\infty$ deformation of the cup product. Equivalently, the cup product is an $A_\infty$ deformation of the $C_\infty$ algebra $C^{\bullet}(A)$ coming from the $\rm{Calc}_\infty $ structure. But, as we discussed above, there may be many such deformations. So far we do not even know whether the DGA $(C^{\bullet}(A), \delta, \smile)$ is $A_\infty$ equivalent to a $C_\infty$ algebra.

 {\bf Noncommutative version of I.} Here the situation becomes much easier: a noncommutative version of I is explicit and canonical. Namely, the negative cyclic complex ${\rm{CC}}^-_{-\bullet}(A)=(C_{-\bullet}(A)[[u]],\,b+uB)$ is an $L_\infty$ module over the DGLA $(\g_A[\epsilon, u], \delta+u\frac{\partial}{\partial \epsilon})$ where $\g_A=(C^{\bullet+1}_A, \delta, [,\,])$ is the DG Lie algebra of Hochschild cochains with the Gerstenhaber bracket. We give explicit formulas for this $L_\infty$ structure (Theorem \ref{thm:1}) and prove that this structure is $L_\infty$ equivalent to the one induced by the ${\rm{Calc}}_\infty$ structure on $(C^\bullet(A), C_\bullet(A))$ (Theorem \ref{thm:2}).

 The reason for the latter statement is the following. Unlike the associative multiplication, the Lie algebra structures that are part of the definition of a calculus live, degreewise, on the very edge of the calculus structure and do not have any room for change. For example, unlike the product, the Lie bracket on a Gerstenhaber algebra cannot be universally deformed, simply because there are no universal operations of needed degrees. Likewise, the DG Lie algebra structure of Ib) cannot be universally deformed: universal operations of needed degrees are too few and easily controlled.

 To define the $L_\infty$ structure from Theorem \ref{thm:1}, we, following \cite{T}, construct it from the noncommutative version of the ring IId) of differential operators on forms, namely the $A_\infty$ algebra ${\rm{CC}}^-_{-\bullet}(C^\bullet(A))$ that was studied in \cite{TT}. It would be interesting, instead of describing the latter by explicit formulas, to interpret it as a part of the $A_\infty$ category of $A_\infty$ functors as described in \cite{KS}, \cite{Kel}, as well as in \cite{Ta}. (In a crude form, this was basically the idea in \cite{NT}).

We conclude the paper by showing that, given a family $\mathcal{A}$ of algebras on a variety $S$, assuming that this family admits a connection as a family of vector spaces, there exists a flat superconnection on the family of complexes $s\mapsto {\rm{CC}}^{\rm{per}}_{-\bullet}({\mathcal{A}}_s).$ This generalizes Getzler's construction of the Gauss-Manin connection \cite{G} from the level of homology to the level of actual complexes. The general fact that the existence of a flat superconnection follows from an existence on an $L_\infty$ module structure is due to Barannikov \cite{Bar}, Remarks 3.3 and 6.7. Note that a modified version of Getzler's construction was used in \cite{KKP}.

 It would be interesting to compare the results of this paper to the ones from \cite{Gin} and \cite{ToV}.

 {\bf Concluding remarks.} We see that noncommutative differential calculus has two levels. At the higher level, there is an inexplicit, noncanonical structure of a ${\rm{Calc}}_\infty$ algebra on the pair $(C^\bullet(A), C_\bullet(A)).$  At the lower level, there are some simpler structures whose existence is implied by the existence of the ${\rm{Calc}}_\infty$ structure but they themselves are explicit and canonical. They are:

 (1) the $A_\infty$ algebra ${\rm{CC}}^-_{-\bullet}(C^\bullet(A))$ and the $A_\infty$ module ${\rm{CC}}^-_{-\bullet}(A)$ over it (noncommutative differential operators on forms);

(2) the DGLA ${\mathfrak g}_A^\bullet = C^{\bullet+1}(A)$ (noncommutative multivectors);

(3) the $L_\infty$ module structure on ${\rm{CC}}^-_{-\bullet}(A)$ over $({\mathfrak g}_A[\epsilon, u], \delta+u\frac{\partial}{\partial \epsilon})$ (noncommutative analog of multivectors acting on forms by Lie derivative and by contraction); and also

(4) the calculus $(H^\bullet(A), H_\bullet(A))$, the homology of $(C^\bullet(A), C_\bullet(A).$ The explicit formulas were given in \cite{DGT}.

A large group of symmetries acts on the space of choices for the ${\rm{Calc}}_\infty$ structure. Note that the group acts not by automorphisms of any structure but on the space of choices of the structure itself. There is an important case when (some extension of) this group acts by automorphisms of the structures (1), (2), (3), and (4).

%automorphisms of any of the structures (1), (2), (3), and (4).
Indeed, take for $A$ the sheaf of functions on a smooth manifold (real, complex, or algebraic). A formality theorem is true for the ${\rm{Calc}}_\infty$ structure above (cf. \cite{DTT2}, \cite{DTT3}); namely, for any choice $\alpha$ of the ${\rm{Calc}}_\infty$ structure, there is a ${\rm{Calc}}_\infty$
quasi-isomorphism of sheaves of ${\rm{Calc}}_\infty$ algebras
\begin{equation}\label{eq:fty}
\Phi _\alpha: (C^\bullet({{\mathcal O}_X}), C_\bullet({{\mathcal O}_X}))_\alpha \isomoto (\wedge^\bullet T_X, \Omega^\bullet _X)
\end{equation}
Here the left hand side is equipped with the ${\rm{Calc}}_\infty$ structure given by $\alpha$ and the right hand side with the standard calculus structure.
%Also, because any of the structures (1)-(4) is canonical up to (quasi-)isomorphism, we have
%\begin{equation}\label{eq:fty 1}
%\Psi _\alpha: (C^\bullet({{\mathcal O}_X}), C_\bullet({{\mathcal O}_X}))_\alpha \isomoto (C^\bullet({{\mathcal O}_X}), C_\bullet({{\mathcal O}_X}))_{\rm{standard}
%\end{equation}
%Then $\Phi_{\alpha\beta}=\Phi_\beta\Psi_\beta^{-1}\Psi_\alpha \Phi_{\alpha}^{-1}$
The cohomology of the right hand side can be identified with the standard calculus
$(H^\bullet({\mathcal O}_X), H_\bullet ({\mathcal O}_X))$; comparing two such identifications,
for any $\alpha$ and $\beta$ we get an automorphism of the standard calculus
\begin{equation}\label{eq:symmetry}
\Phi_{\alpha \beta}:(H^\bullet (X, \wedge^\bullet T_X); H^\bullet(X, \Omega^\bullet _X))\isomoto (H^\bullet (X, \wedge^\bullet T_X); H^\bullet(X, \Omega^\bullet _X))
\end{equation}

It would be very interesting to compare the above construction to the $L_\infty$ quasi-isomorphisms constructed by Merkulov in \cite{Me1}.

In \cite{K1}, Kontsevich constructed automorphisms of the cohomology with coefficients
in multivector fields that are probably
part of the above construction. Similarly, for any $\alpha$ and $\beta$
one constructs an automorphism of the classical analog of any of the structures (1)-(3): an $L_\infty$ quasi-isomorphism of $\Omega ^{0,\bullet} (X, \wedge ^{\bullet+1} T_X)$ with itself; a compatible quasi-isomorphism of $L_\infty$ modules over $(\Omega ^{0,\bullet} (X, \wedge ^{\bullet+1} T_X)[\epsilon, u], {\overline{\partial}}+u\frac{\partial}{\partial\epsilon})$ between $(\Omega ^{0,\bullet} (X,\Omega ^{-\bullet}_X)[[u]], {\overline{\partial}}+u\partial)$ and itself, etc. Indeed, one has
$$C^{\bullet+1}({\mathcal O}_X)_{\rm{Gerst}} \stackrel{\Psi_\alpha}{\longleftarrow} C^{\bullet+1}({\mathcal O}_X)_\alpha \stackrel{\Phi_\alpha}{\longrightarrow} \wedge^{\bullet+1} T_X$$
The DGLA on the left is the Hochschild cochain complex with the Gerstenhaber bracket. The map on the right is a $G_\infty$ quasi-isomorphism of sheaves of $G_\infty$ algebras (formality); the one on the right is the $L_\infty$ quasi-isomorphism of sheaves of DGLA (rigidity). Comparing two such sequences for $\alpha$ and $\beta$, we get an $L_\infty$ quasi-isomorphism at the level of algebras; similarly for modules.

Let us finish by some remarks about the algebraic index theorem. The formality ${\rm{Calc}}_\infty$ quasi-isomorphism implies a quasi-isomorphism of complexes
\begin{equation}\label{eq:index}
\Omega ^{0,\bullet} (X,{\rm{CC}}^{\rm{per}}_{-\bullet}({\mathcal O}_X))\isomoto \Omega ^{0,\bullet} (X,\Omega_X^{-\bullet}((u)))
\end{equation}
An algebraic index theorem is a statement comparing it to the standard Hochschild-Kostant-Rosenberg map (and an analogous statement for a deformation quantization of ${\mathcal O}_X$; cf. \cite{BNT}, \cite{bgnt1} for the symplectic case). Equivalently, it is a statement about the image of the zero-homology class $1$. One can show that this image is an expression in the Chern classes of $T_X$ that becomes $\sqrt{{\widehat{A}}(T_X)}$ if we send all the odd Chern classes to zero. If the automorphisms in \cite{K1} do extend to automorphisms of calculi that come from \eqref{eq:symmetry}, then it looks like any multiplicative characteristic class with the above property may occur, for an appropriate choice of $\alpha.$ Note also that nothing in our argument implies that the symmetries constructed above are   ${\rm{Calc}}_\infty$ quasi-isomorphisms of $(\wedge ^\bullet T_X, \Omega ^\bullet _X)$ with itself. On the other hand, there definitely are some ${\rm{Calc}}_\infty$ quasi-isomorphisms, for instance the exponential of the following derivation:
$$\iota _{c_1(TX)}:\Omega ^{0,\bullet} (X,\wedge ^{\bullet}T_X)\to \Omega ^{0,\bullet+1} (X,\wedge ^{\bullet-1}T_X)$$
and
$${c_1(TX)}\wedge:\Omega ^{0,\bullet} (X,\Omega^\bullet_X)\to \Omega ^{0,\bullet+1} (X,\Omega^{\bullet+1}_X)$$

Note also that the Hochschild-Rosenberg map, followed with multiplication by the square root of the Todd class, appears in \cite{Ma}, \cite{Ca1}, \cite{Ca2}, \cite{Ca3}, and \cite{Ra} and is characterized by preserving another algebraic structure. Namely, it intertwines the Mukai pairing with teh standard pairing on the cohomology. This suggests that the correct formality theorem could be formulated for an algebraic structure encompassing both the calculus and the pairing, probably related to (genus zero part of) the TQFT structure from \cite{Co} and \cite{KS}.

{\bf Acknowledgments.}
D.T. and B.T. are supported by
NSF grants. The work of V.D. is partially supported by the Grant
for Support of Scientific Schools NSh-3036.2008.2.
We are grateful to Paul Bressler, Kevin Costello, Ezra Getzler, Maxim Kontsevich, and Yan Soibelman for fruitful discussions.

 \section{Operators on forms in noncommutative calculus}\label{s:oofincc}
 \subsection{The Hochschild cochain complex} \label{hocochain}

Let $A$ be a graded algebra with unit over a commutative unital ring
$K$ of characteristic zero.  A Hochschild $d$-cochain is a linear map $A^{\otimes d}\to A$.  Put,
for $d\geq 0$,
$$
        C^d(A) = C^d (A,A) =\operatorname{Hom}_K({\overline{A}}^{\otimes d},A)
$$
where ${\overline{A}}=A/K\cdot 1$. Put
$|D|=($ degree of the linear map $D)+d$

  Put for cochains $D$ and $E$ from $C^{\bullet}(A,A)$
$$
        (D\smile E)(a_1,\dots,a_{d+e})=(-1)^{| E|\sum_{i \leq d}(|a_i| + 1)}
        D(a_1,\dots,a_d) E(a_{d+1},\dots,a_{d+e});
$$
$$
        (D\circ E)(a_1,\dots,a_{d+e-1})=\sum_{j \geq 0}
        (-1)^{(|E|+1)\sum_{i=1}^{j}(|a_i|+1)}
D(a_1,\dots,a_j,
        E(a_{j+1},\dots,a_{j+e}),\dots);
$$
$$
        [D, \; E]= D\circ E - (-1)^{(|D|+1)(|E|+1)}E\circ D
$$
These operations define the graded associative algebra
$(C^{\bullet}(A,A)\;,\smile)$ and the graded Lie algebra
($C^{\bullet + 1}(A,A)$, $[\;,\;]$) (cf. \cite{CE}; \cite{G}).
Let
$$
        m(a_1,a_2)=(-1)^{\deg a_1}\;a_1 a_2;
$$
this is a 2-cochain of $A$ (not in $C^2$).  Put
$$
        \delta D=[m,D];
$$
$$
        (\delta D)(a_1,\dots,a_{d+1})=(-1)^{|a_1||D|+|D|+1}
 a_1 D(a_2,\dots,a_{d+1})+
$$
$$
        +\sum
_{j=1}^{d}(-1)^{|D|+1+\sum_{i=1}^{j}(|a_i|+1)}
                D(a_1,\dots,a_ja_{j+1},\dots,a_{d+1})
 +(-1)^{|D|\sum_{i=1}^{d }(|a_i|+1)}D(a_1,\dots,a_d)a_{d+1}
$$

One has
$$
        \delta^2=0;\quad\delta(D\smile E)=\delta D\smile E+(-1)^{|deg D|}
                D\smile\delta E
$$
$$
        \delta[D,E]=[\delta D,E]+(-1)^{|D|+1}\;[D,\delta E]
$$
($\delta^2=0$ follows from $[m,m]=0$).

Thus $ C^{\bullet}(A,A)$ becomes a complex; we will denote it also by $C^{\bullet}(A)$. The cohomology of this complex
is $H^{\bullet}(A,A)$ or the Hochschild cohomology. We denote it also by
$H^{\bullet}(A) $.  The $\smile$ product induces the
Yoneda product on $H^{\bullet}(A,A)=Ext_{A\otimes A^0}^{\bullet}(A,A)$.  The operation
$[\;,\;]$ is the Gerstenhaber bracket \cite{Ge}.

If $(A, \;\; \partial)$ is a differential graded algebra then one can define
the differential $\partial$ acting on $A$ by
$$
\partial D \;\; = \; [\partial , D]
$$

\begin{teo} \cite{Ge} The cup product and the Gerstenhaber bracket induce a Gerstenhaber algebra structure on $H^\bullet(A)$.
\end{teo}
%%%%%%%%%%%%%%%%%%%%%%%%%%%%%%%%%5%%%%%%%%%%%%%%%%%%%%%%%%%%%%%%5%%%%%%%%%%%%%%%%%%%%%%%%%%%%%%5%%%%%%%%%%%%%%%%%%%%%%%%%%%%%%5%%%%%%%%%%%%%%%%%%%%%%%%%%%%%%5%%%%%%%%%%%%%%%%%%%%%%%%%%%%%%5%%%%%%%%%%%%%%%%%%%%%%%%%%%%%%5%%%%%%%%%%%%%%%%%%%%%%%%%%%%%%5%%%%%%%%%%%%%%%%%%%%%%%%%%%%%%5%%%%%%%%%%%%%%%%%%%%%%%%%%%%%%5%%%%%%%%%%%%%%%%%%%%%%%%%%%%%%5%%%%%%%%%%%%%%%%%%%%%%%%%%%%%%5%%%%%%%%%%%%%%%%%%%%%%%%%%%%%%5%%%%%%%%
\subsection{Hochschild chains}  \label{ss:hochschild-1}
%%%%%%%%%%%%%5%%%%%%%%%%%%%%%%%%%%%%%%%%%%%%5%%%%%%%%%%%%%%%%%%%%%%%%%%%%%%5%%%%%%%%%%%%%%%%%%%%%%%%%%%%%%5%%%%%%%%%%%%%%%%%%%%%%%%%%%%%%5%%%%%%%%%%%%%%%%%%%%%%%%%%%%%%5%%%%%%%%%%%%%%%%%%%%%%%%%%%%%%5%%%%%%%%%%%%%%%%%%%%%%%%%%%%%%5%%%%%%%%%%%%%%%%%%%%%%%%%%%%%%5%%%%%%%%%%%%%%%%%%%%%%%%%%%%%%5%%%%%%%%%%%%%%%%%%%%%%%%%%%%%%5%%%%%%%%%%%%%%%%%%%%%%%%%%%%%%5%%%%%%%%%%%%%%%%%%%%%%%%%%%%%%5%%%%%%%%%%%%%%%%%%%%%%%%%%%%%%5%%%%%%%%%%%%%%%%%%%%%%%%%%%%%%%5%%%%%%%%%%%%%%%%%%%%%%%%%%5
%%%%%%%%%%%%%%%%%%%%%%%%%%%%%%5
Let $A$ be an associative unital DG algebra over a ground ring $K$. The differential on $A$ is denoted by $\delta$. Recall that by definition
$$\overline{A} = A / K\cdot 1$$
Set
$$C_p (A,A) = C_p(A) = A \otimes \overline{A} ^{\otimes p}$$
Define the differentials $\delta: C_{\bullet}(A) \to C_{\bullet}(A)$, $b: C_{\bullet}(A) \to C_{\bullet - 1}(A)$, $B: C_{\bullet}(A) \to C_{\bullet + 1}(A)$ as follows.
\[
\delta (a_0\otimes\cdots\otimes a_p ) =
\sum_{i=1}^p {(-1)^{\sum_{k<i}{(| a_k| + 1)+1}}
(a_0\otimes\cdots\otimes\delta a_i \otimes \cdots \otimes a_p )};
\]

\begin{equation} \label{eq:b grad}
b(a_0 \otimes \ldots \otimes a_p) = \sum _{k=0}^{p-1} (-1)^{\sum_{i=0}^{k} {(|a_i| + 1)+1}}
 a_0 \ldots \otimes a_k a_{k+1} \otimes \ldots a_p
\end{equation}
$$+ (-1)^{|a_p| + (|a_p|+1)\sum_{i=0}^{p-1}(|a_i|+1)} a_pa_0 \otimes \ldots \otimes a_{p-1};
$$

\begin{equation} \label{eq: B graded}
B(a_0 \otimes \ldots \otimes a_p) = \sum_{k=0}^p (-1)^{\sum _{i \leq k}(|a_i| + 1) \sum _{i \geq k}(|a_i| + 1)} 1 \otimes a_{k+1} \otimes \ldots a_p \otimes
a_0 \otimes \ldots \otimes a_k
\end{equation}
The complex $C_{\bullet}(A)$ is the total complex of
the double complex with the differential $b + \delta$.

Let $u$ be a formal variable of degree two. The complex $(\Ccc [[u]], b+\delta + uB)$ is called {\it {the negative cyclic complex}} of $A$.

Now put
\begin{equation} \label{eq: L}
L_D(a_0 \otimes \ldots \otimes a_n)=\sum _{k=1}^{n-d} \epsilon _k a_0 \otimes \ldots \otimes D(a_{k+1}, \ldots, a_{k+d}) \otimes \ldots \otimes a_n +
\end{equation}
$$ \sum _{k=n+1 -d}^{n} \eta _k D (a_{k+1}, \ldots, a_n, a_0, \ldots ) \otimes \ldots \otimes a_k
$$

(The second sum in the above formula is taken over all cyclic permutations such that $a_0$ is inside $D$). The signs are given by
$$ \epsilon _k = (|D| + 1)(|a_0|+\sum _{i=1}^{k} (|a_i| +1))$$
and
$$
\eta _k = |D|+ \sum_{i \leq k}(|a_i|+1)\sum_{i \geq k}(|a_i|+1)
$$
\begin{proposition}
$$[L_D, L_E]=L_{[D,E]};\;\;\;
[b, L_D] + L_{\delta D} = 0;\;\;\;
[L_D, B] = 0$$
\end{proposition}

 \subsection{$A_\infty$ algebras and modules} \label{ss:aialia}
 Recall \cite{LS}, \cite{St} that an $A_{\infty}$ algebra is a graded vector space ${\cal{C}}$ together with a Hochschild cochain $m$ of total degree $1$,
$$ m = \sum _{n=1}^{\infty} m_n$$
where $m_n \in C^n({\cal{C}})$ and
$$[m,m] = 0$$
Recall also the definition of $A_{\infty}$ modules over $A_{\infty}$ algebras. First, note that for a graded space ${\cal{M}}$, the Gerstenhaber bracket $[\;,\;]$ can be extended to the space
$$\operatorname{Hom} (\overline {{\cal{C}}}^{\otimes {\bullet}}, {\cal{C}}) \oplus \operatorname{Hom} ({\cal{M}} \otimes \overline {{\cal{C}}}^{\otimes {\bullet}}, {\cal{M}})
$$

For a graded $k$-module ${\cal{M}}$, a structure of an $A_{\infty}$ module over an $A_{\infty}$ algebra ${\cal{C}}$ on ${\cal{M}}$ is a cochain of total degree one
$$ \mu = \sum _{n=1}^{\infty} \mu _n$$
$$ \mu _ n \in {Hom} ({\cal{M}} \otimes \overline {{\cal{C}}}^{\otimes {n-1}}, {\cal{M}})$$
such that
$$[m+\mu, m+\mu ] =0$$

 \subsection{The $A_\infty$ structure on chains of cochains}\label{ss:taisococ}
\begin{thm} \label{thm: a-infty mod}
There is a structure $\{{\bf m}_n\}$ of an $A_\infty$ algebra on $CC^-_{-\bullet}(C^\bullet(A)),$ and a structure $\{\mu_n\}$ of an $A_{\infty}$ module over this $A_{\infty}$ algebra on $C_{-\bullet}(A)[[u]],$ such that:
\begin{itemize}
\item All ${\bf m}_n$ and $\mu_n$ are $k[[u]]$-linear, $(u)$-adically continuous.
\item ${\bf m}_1=b+\delta+uB;$ $\mu_1 = b +uB.$
\item Modulo $u$, the space $C_0(C^\bullet(A))=C^\bullet(A)$ is a subalgebra, with the structure given by
 the cup product.
%\item $\mu_2 (a, D) = (-1)^{|a||D| + |a|}(i_D + u S_D)a$
\item For $a \in C_{-\bullet}(A)[[u]],$ $D\in C^\bullet(A)$:$\;\;\mu_2 (a, 1 \otimes D) = (-1)^{|a||D|}L_D a. $

%For $a, \; x \in C_{-\bullet}(A)[[u]]$:
%$(-1)^{|a|}\mu_2(a,x) = (\operatorname{sh} + u \operatorname{sh}')(a,x)$
\end{itemize}
\end{thm}

 Explicit formulas can be found in \cite{TT}, \cite{T}. The proof is given in \cite{T}.

 \section{The $L_\infty$ module structure on the negative cyclic complex}\label{s:tlimsotncc}
 Now introduce the following differential graded algebras. Let $C(\gA[u, \epsilon])$ be the standard Chevalley-Eilenberg chain complex of the DGLA $\gA[u, \epsilon]$ over the ring of scalars $K[u]$. It carries the Chevalley-Eilenberg differential $\partial$ and the differentials $\delta$ and $\partial _{\epsilon}$ induced by the corresponding differentials on $\gA[u, \epsilon]$. Let $C_+(\gA[u, \epsilon])$ be the augmentation co-ideal, i.e. the sum of all positive exterior powers of our DGLA. The comultiplication defines maps
$$
C_+(\gA[u, \epsilon])\mapsto C_+(\gA[u, \epsilon])^{\otimes n};
$$
$$
c\mapsto \sum c_1^+\otimes \ldots \otimes c_n^+
$$
\begin{definition}\label{dfn:Bar}
Define the associative DGA $B(\gA[u, \epsilon])$ over $K[[u]]$ as the tensor algebra of $C_+(\gA[u, \epsilon])$ with the differential $d$ determined by
$$dc = (\delta + \partial)c-\frac{1}{2}\sum (-1)^{|c_1^+|}c_1^+c_2^+ + u\partial _{\epsilon}c.$$
\end{definition}
\begin{definition}\label{dfn:Bar tw}
Let the associative DGA $B^{\operatorname {tw}}(\gA[u, \epsilon])$ over $K[[u]]$ be the tensor algebra of $C_+(\gA[u, \epsilon])$ with the differential $d$ determined by
$$dc = (\delta + \partial)c -\frac{1}{2}\sum (-1)^{|c_1^+|}c_1^+c_2^+ +u\sum_{n=1}^{\infty} \partial _{\epsilon}c_1^+\ldots \partial _{\epsilon}c_n^+.$$
\end{definition}
A structure of an $L_\infty$ module over $\gA[u, \epsilon]$ on a complex ${\cal{M}}$ is by definition a morphism of DGA $B(\gA[u, \epsilon])\to{\rm{End}}({\cal{M}}).$ It would be nice to have explicit formulas for such a morphism. What we can do instead is construct an explicit morphism
$$B^{\operatorname {tw}}(\gA[u, \epsilon])\to{\rm{End}}{\rm{CC}}_{-\bullet}^-(A)$$
together with a quasi-isomorphism of DGAs
$$U(\gA[u, \epsilon])\rightarrow B^{\operatorname {tw}}(\gA[u, \epsilon]).$$
%We start with a key property of the $\Ai$ structures from section \ref{ss:taisococ}.
%\begin{proposition}\label{prop:key property 1}
%Both ${\bf m}_k(c_1, \ldots, c_k)$ and $\mu _k(c_1, \ldots, c_k)$ are equal to zero if one of the arguments $c_i$, $i<k$, is of the form $1\otimes \ldots$.
%\end{proposition}
%This can be seen immediately from the definition \cite{TT}, \cite{T}.

Let $S(\gA)^+$ be the augmentation ideal, and let
\begin{equation}\label{eq:coproduct}
Y\mapsto \sum Y^+_1\otimes \ldots \otimes Y_n^+
\end{equation}
denote the map
\begin{equation}\label{eq:coproduct 1}
S(\gA)^+ \to (S(\gA)^+)^{\otimes n}
\end{equation}
defined as the n-fold coproduct, followed by the $n$th power of the projection from $S(\gA)$ to $S(\gA)^+$ along $K\cdot 1$.
\begin{definition}\label{dfn:action of btw}
For $n\geq 1$, define:
$$x\cdot(\epsilon E_1 \wedge \ldots \wedge\epsilon E_n)=\sum_{n\geq 1}(-1)^{|x|}\mu _{n+1}(x, \, {\overline{Y}}_1^+, \ldots , {\overline{Y}}_n^+); $$
for $n\geq 0$,
$$x\cdot(\epsilon E_1 \wedge \ldots \wedge\epsilon E_n \wedge D)= \sum_{n\geq 1}(-1)^{|x|}\mu _{n+2}(x, \, {\overline{Y}}_1^+, \ldots , {\overline{Y}}_n^+, 1\otimes D);$$
$$x\cdot(  \epsilon E_1 \wedge \ldots \wedge\epsilon E_n\wedge D_1 \wedge \ldots \wedge D_k)= 0$$ for $k>1.$
Here $D, \; D_i, \; E_j \in \gA$ and $Y=E_1 \ldots E_n \in S(\gA)^+.$
\end{definition}
\begin{pred}\label{pred:action of btw}
The formulas from Definition \ref{dfn:action of btw} above define an action of the DGA $B^{\operatorname {tw}}(\gA[u, \epsilon])$ on ${\rm{CC}}^-_{-\bullet}(A).$
\end{pred}
The proof is contained in \cite{T}.
%%%%%%%%%%%%%%%%%%%%%%%%%%%%%%%%%%%%%%%%%%%%%%%%%%%%%%%%%%%%%%%%%%%%%%%%%%%%%%%%%%%%%%%%%%%%%%%%%%%%%%%%%%%%%%%%%%%%%%%%%%%%%%%%%%%%%%%%%%%%%%%%%%%%%%%%%%%%%%%%%%%%%%%%%%%%%%%%%%%%%%%%%%%%%%%%%%%%%%%%%%%%%%%%%%%%%%%%%%%%%%%%%%%%%%%%%%%%%%%%%%%%%%%%%%%%%%%%%%%%%%%%%%%%%%%%%%%%%%%%%%%%%%%%%%%%%%%%%%%%%%%%%%%%%%%%%%%%%%%%%%%%%%%%%%%%%%%%%%%%%%%%%%%%%%%%%%%%%%%%%%%%%%%%%%%%%%%%%%
\subsection{The $L_\infty$ action}\label{tlia}
%%%%%%%%%%%%%%%%%%%%%%%%%%%%%%%%%%%%%%%%%%%%%%%%%%%%%%%%%%%%%%%%%%%%%%%%%%%%%%%%%%%%%%%%%%%%%%%%%%%%%%%%%%%%%%%%%%%%%%%%%%%%%%%%%%%%%%%%%%%%%%%%%%%%%%%%%%%%%%%%%%%%%%%%%%%%%%%%%%%%%%%%%%%%%%%%%%%%%%%%%%%%%%%%%%%%%%%%%%%%%%%%%%%%%%%%%%%%%%%%%%%%%%%%%%%%%%%%%%%%%%%%%%%%%%%%%%%%%%%%%%%%%%%%%%%%%%%%%%%%%%%%%%%%%%%%%%%%%%%%%%%%%%%%%%%%%%%%%%%%%%%%%%%%%%%%%%%%%%%%%%%%%%%%%%%%%%%%%%
It remains to pass from $B^{\operatorname {tw}}(\gA[\epsilon, u])$ to $U(\gA[\epsilon, u])$.

The following is contained in \cite{T}.
\begin{lemma}\label{quis of twisted and untwisted}
The formulas
$$D\to D;$$
$$\epsilon E_1 \wedge \ldots \wedge \epsilon E_n \mapsto \frac{1}{n!}\sum_{\sigma \in S_n}\frac{1}{n!}(\epsilon E_{\sigma_1}) E_{\sigma_2}\ldots E_{\sigma_n};$$
$$D_1\wedge \ldots D_k \wedge \epsilon E_1 \wedge \ldots \wedge \epsilon E_n\mapsto 0$$
for $k>1$ or $k=1, \, n\geq 1$
define a quasi-isomorphism of DGAs
$$B^{\operatorname {tw}}(\gA[\epsilon, u]) \to U(\gA[\epsilon, u]).$$
\end{lemma}
{\bf{Proof.}} The fact that the above map is a morphism of DGAs follows from an easy direct computation.
To show that this is a quasi-isomorphism, consider the increasing filtration by powers of $\epsilon$. At the level of graduate quotients, $B^{\operatorname {tw}}(\gA[\epsilon, u])$ becomes the standard free resolution of $(U(\gA[\epsilon, u]), \delta)$, and the morphism is the standard map from the resolution to the algebra, therefore a quasi-isomorphism. The statement now follows from the comparison argument for spectral sequences.

To summarize, we have constructed explicitly a DGA $B^{\operatorname {tw}}(\gA[\epsilon, u])$ and the morphisms of DGAs
\begin{equation}\label{eq:thm2} U(\gA[\epsilon, u]) \leftarrow B^{\operatorname {tw}}(\gA[\epsilon, u]) \to \operatorname{End} _{K[[u]]}(\CCnA)
\end{equation}
where the morphism on the left is a quasi-isomorphism. This yields an $A_{\infty}$ morphism
$$U(\gA[\epsilon, u])\to \operatorname{End} _{K[[u]]}(\CCnA)$$
and therefore an $L_{\infty}$ morphism
$$\gA[\epsilon, u] \to \operatorname{End} _{K[[u]]}(\CCnA).$$
We get
\begin{teo}\label{thm:1} The maps (\ref{eq:thm2})
define on ${\rm{CC}}^-_{-\bullet}(A)$ a $K[u]$-linear, $(u)$-adically continuous structure of an $L_\infty$ module over the DGLA $(\g _A[\epsilon, u], \delta+u\frac{\partial}{\partial \epsilon})$ where $\g _A=(C^{\bullet+1}(A), \delta, [,]_G).$
 \end{teo}
 \begin{rmk} The property of being $K[u]$-linear is crucial. Indeed, as pointed out by K. Costello, $\g[\epsilon, u]$ is quasi-isomorphic to $\g,$ so every $\g$-module is an $L_\infty$ module over $\g[\epsilon, u]$ by transfer of structure.
 \end{rmk}
 
 In Section \ref{s:tgmc} we will need the following strengthening of Theorem \ref{thm:1}.
\begin{lemma}\label{lemma:relative cochain}
On ${\rm{CC}}^-_{-\bullet} (A)$, there exists a $K[u]$-linear, $(u)$-adically continuous structure of an $L_\infty$ module over the DGLA $(\g _A[\epsilon, u], \delta+u\frac{\partial}{\partial \epsilon})$ where $\g _A=(C^{\bullet+1}(A), \delta, [,]_G)$, such that, for any derivation $D$ of $A$, any $n>1$ and any $X_2, \ldots, X_n$ in $\g _A[\epsilon, u]$, $D\wedge X_2\wedge \ldots \wedge X_n \mapsto 0.$
\end{lemma}

{\bf Proof.} The statement is almost obvious. Indeed, all the above constructions are clearly invariant under the Lie algebra of infinitesimal symmetries ${\rm{Der}}(A)$, and the image of $D\wedge X_2\wedge \ldots \wedge X_n$ is the homotopy for the expression that measures failure of the previous maps to be invariant. The actual proof goes as follows. Start by observing that the negative cyclic complex carries an action of a DGA smaller than $U(B(\g _A[\epsilon, u]))$ or $U(B^{\rm{tw}}(\g _A[\epsilon, u]))$. Namely, for any DGLA $\g$, denote by $I$ the DG ideal ${\rm{Ker}}(U(B(\g))\to U(\g)).$ Note that $B(\g)$ is a DG subalgebra both of $B(\g [\epsilon, u])$ and of $B^{\rm{tw}}(\g [\epsilon, u])$. Denote by ${\overline U}(B(\g [\epsilon, u]))$, resp. by  ${\overline U}(B^{\rm{tw}}(\g [\epsilon, u]))$, the quotient of the corresponding DGA by the two-sided ideal generated by $I$. Both are DGAs that are free over $U(\g)$ as graded algebras. Moreover, they are $U(\g)$-semifree resolutions of $U(\g[\epsilon, u]).$ It is clear from the constructions that ${\overline U}(B(\g_A [\epsilon, u]))$ and ${\overline U}(B^{\rm{tw}}(\g_A [\epsilon, u]))$ act on the negative cyclic complex through their quotients ${\overline U}.$ Indeed, the image of $I$ in $U(B^{\rm{tw}})$ acts by zero because $D_1\wedge \ldots \wedge D_n$ do, $D_i$ in $\g_A$, $n>1;$ and there is a morphism
\begin{equation}\label{eq:morphism red}
  {\overline U}(B(\g [\epsilon, u]))\to {\overline U}(B^{\rm{tw}}(\g [\epsilon, u]))
 \end{equation}
 over $U(\g[\epsilon, u])$, since both are $U(\g)$-semifree resolutions.

Denote by $L_D$ the natural action of $\g$ on ${\overline U}(B^{\rm{tw}}(\g [\epsilon, u]))$ and ${\overline U}(B^(\g [\epsilon, u]))$:
$$L_D(X_1\wedge \ldots \wedge X_n)=[D, X_1\wedge \ldots \wedge X_n]-\sum_k \pm  (X_1\wedge \ldots [D, X_k]\ldots\wedge X_n).$$
Here the first commutator is taken in the associative algebra ${\overline U}.$ Define also
$$\iota_D(X_1\wedge \ldots \wedge X_n)=D\wedge X_1\wedge \ldots \wedge X_n.$$
The above operators satisfy the usual relations
$$[L_{D}, L_{E}]=L_{[D,E]};\;[L_{D}, \iota_{E}]=\iota_{[D,E]};\;[\iota_{D}, \iota_{E}]=0;\;[d,\iota_D]=L_D$$
where $d$ is the total differential in ${\overline U}(B^{\rm{tw}}(\g [\epsilon, u]))$, resp. in ${\overline U}(B(\g [\epsilon, u]))$. It is clear that the morphism (\ref{eq:morphism red}), being given by universal formulas in terms of the commutator, is invariant with respect to $L_D.$ In fact it can be chosen to be invariant also with respect to $\iota_D.$ This is easy to see using the standard inductive procedure of defining this morphism on free generators and the fact that the space of generators is a free module over the free commutative algebra generated by all $\iota_D.$

Finally, we conclude that the morphism (\ref{eq:morphism red}) sends the ideal generated by $D\wedge X_2\wedge\ldots X_n$ as in the statement of the Lemma to the analogous ideal ${\overline U}(B^{\rm{tw}}).$ It remains to show that the latter acts on the negative cyclic complex by zero. This follows from the second formula of Definition \ref{dfn:action of btw} and from the fact that, in the notation of Theorem \ref{thm: a-infty mod}, $\mu_n(x, a_1, \ldots, a_n, 1\otimes D)=0$ for $D\in C^1(A,A)$ and $n\geq 1.$ The latter follows from the explicit formulas in \cite{TT}, \cite{T}. This ends the proof of the Lemma.
 \section{Relation to the homotopy calculus structure}\label{s:rtthcs}
 \subsection{Calculi}\label{ss:c}
 A {\it Gerstenhaber algebra} is a graded space ${\cal{V}} ^{\bullet}$
together with
\begin{itemize}
\item A graded commutative associative algebra structure on
${\cal{V}}^{\bullet}$;
\item a graded Lie algebra structure on ${\cal{V}}^{{\bullet}+1}$
such that
$$[a,bc]=[a,b]c+(-1)^{(|a| -1) |b|)}b[a,c]$$
\end{itemize}
 \begin{definition} \label{dfn:precalc}
A {\it precalculus} is a pair of a Gerstenhaber algebra
${\cal{V}}^{\bullet}$ and a graded space $\Omega ^{\bullet}$  together with
\begin{itemize}
\item a structure of a graded module over the graded commutative  algebra
${\cal{V}}^{\bullet}$ on $\Omega ^{-{\bullet}} $ (corresponding action is
denoted by $i_a,\; a \in {\cal{V}}^{\bullet}$);
\item a structure of a graded module over the graded Lie  algebra
${\cal{V}}^{\bullet +1}$ on ${\Omega} ^{-{\bullet}}$ (corresponding action
is denoted by $L_a,\; a \in {\cal{V}}^{\bullet}$)
such that
$$[i_a,L_b]=i_{[a,b]}$$
and
$$L_{ab} = L_a i_b + (-1)^{|a|} i_a L_b$$
\end{itemize}
\end{definition}
\begin{definition} \label{dfn:calc}
A {\it calculus} is a precalculus
together with an operator $d$ of degree 1 on $\Omega ^{{\bullet}}$ such that
$d^2 = 0$ and
$$ [d,i_a]=L_a. $$
\end{definition}
\begin{example} \label{ex:calc-M}
For any manifold one defines a calculus $\Ca (M)$ with
${\cal{V}}^{\bullet}$ being the algebra of multivector fields,
$\Omega ^{\bullet}$ the space of differential forms, and $d$ the
De Rham differential. The operator $i_a$ is the contraction of a
form by a multivector field.
\end{example}
\begin{definition}\label{dfn:dg calc} A differential graded calculus is a calculus $({\mathcal V}^\bullet , \Omega^{-\bullet})$ with differentials $\delta$ on ${\mathcal V}^\bullet$ and $b$ on $\Omega^{\bullet}$ that are both of degree one and are derivations with respect to the calculus structure.
\end{definition}
The following construction is motivated by Example
\ref{ex:calc-M}. For a Gerstenhaber algebra ${\cal{V}}^{\bullet}$,
let $\Y ({\cal{V}}^{\bullet})$ be the associative algebra
generated by two sets of generators $i_a$, $L_a$, $a \in
{\cal{V}}^{\bullet}$, both $i$ and $L$ linear in $a$,
$$|i_a| = |a|; \; |L_a| = |a| - 1$$
subject to relations
$$i_ai_b = i_{ab};\;\;[L_a,L_b] = L_{[a,b]};$$
$$[i_a, L_b] = i_{[a,b]};\;L_{ab} = L_ai_b + (-1)^{|a|}i_aL_b$$

The algebra $\Y ({\cal{V}}^{\bullet})$ is equipped with the differential $d$
of degree one which is defined as a derivation sending $i_a$ to $L_a$ and
$L_a$ to zero.

For a smooth manifold $M$ one has a homomorphism
$$\Y (\VM) \rightarrow D (\Omega ^{\bullet}(M))$$
The right hand side is the algebra of differential operators on
differential forms on $M$, and the above homomorphism sends the
generators $i_a$, $L_a$ to corresponding differential operators on
forms (cf. Example  \ref{ex:calc-M}). The
above map is in fact an isomorphism, cf. \cite{DTT2}, Proposition 11 in section 6.3.
 \subsection{Comparing the two $L_\infty$ module structures}\label{ss:cttlims}
 For a DG calculus $({\cal{V}}^{\bullet}, \Omega^{\bullet}),$ let $\g$ be the DGLA $({\cal{V}}^{\bullet+1}, \delta, \{,\}).$ One has:

a) $(\Omega^{-\bullet}[[u]], b+ud)$ is a DG module over $\g.$ Moreover:

b) $(\Omega^{-\bullet}[[u]], b+uB)$ is a DG module over $(\g[\epsilon, u], \delta + u\frac{\partial}{\partial \epsilon}),$ $X+\epsilon Y$ acting via $L_X+\iota _Y.$

The same is true for a ${\rm{Calc}}_\infty$ algebra $ ({\cal{V}}^{\bullet}, \Omega^{\bullet})$ if one replaces DG modules by $L_\infty$ modules. Recall from \cite{DTT2}:

\begin{teo}\label{thm:dtt2} $(C^\bullet(A), C_{\bullet}(A))$ is a ${\rm{Calc}}_\infty$ algebra whose underlying $L_\infty$ structure as in a) above is:
$\g = \g_A= (C^{\bullet+1}(A), \delta, [,]_G),$ acting on ${\rm{CC}}^-_{-\bullet}(A)$ via the Lie derivative $L_D.$
\end{teo}
From this, and from b) above, we conclude that ${\rm{CC}}^-_{-\bullet}(A)$ is an $L_\infty$ module over $(\g[\epsilon, u], \delta + u\frac{\partial}{\partial \epsilon}).$ Indeed, as explained in the introduction, the pair $(C^\bullet, C_\bullet)$ is quasi-isomorphic to another pair of complexes that is actually a DG calculus, with the DGLA structure equivalent to the one given by the Gerstenhaber bracket. Apply b) to that pair and then get the $L_\infty$ module structure on $C_\bullet(A)[[u]]$ by transfer of structure. Note also that an $L_\infty$ algebra and an $L_\infty$ module are in particular complexes, and the differentials coincide with the Hochschild differentials $\delta$ and $b$.
\begin{teo}\label{thm:2}
The above $L_\infty$ structure is equivalent to the one given by Theorem \ref{thm:1}.
\end{teo}
{\bf Proof}. First, recall from \cite{KS}, \cite{KS1}, \cite{TT} the notion of a two-colored operad and a chain of quasi-isomorphisms of two-colored operads
\begin{equation}\label{eq:quises of 2-col ops}
{\rm{Calc}}_{\rm{alg}}\leftarrow {\rm{Calc}}_{\rm{geom}}\leftarrow {\rm{Calc}}_{\infty}\rightarrow {\rm{Calc}}
\end{equation}
Here $Calc_{alg}$ is the operad which
acts on cochains and chain
and which is generated by the cup-product
on cochains, the insertions of cochains into a
cochain, and insertions of components of a chain into
a cochain compatible with the cyclic order on these components.
The precise description of this operad can be found in
\cite{DTT2} (Sect. 4.1) or \cite{KS} (Sect.  11.1, 11.2
and 11.3). Note that $Calc_{alg}$ was denoted by ${\rm{KS}}$ in \cite{DTT2}.

The two-colored operad ${\rm{Calc}}_{\rm{geom}}$ is the operads of chain complexes of the spaces of configurations of little discs (on a disc and on a cylinder); algebras over ${\rm{Calc}}$ are by definition calculi, and ${\rm{Calc}}_{\infty}$ is a cofibrant resolution of ${\rm{Calc}}$; algebras over it are ${\rm{Calc}}_{\infty}$-algebras. The two-colored operad ${\rm{Lie}}_\infty^+$ maps to ${\rm{Calc}}_{\infty}$. An algebra over ${\rm{Lie}}_\infty^+$ is a pair consisting of an $L_\infty$ algebra and an $L_\infty$ module over it. There is a parallel diagram for precalculi.

Note that for any $L_\infty$ algebra ${\mathcal L}$ and any graded space ${\mathcal M}$, an $L_\infty$ ${\mathcal L}$-module structure on ${\mathcal M}$ is the same as a Maurer-Cartan element of the Chevalley-Eilenberg complex $C^\bullet({\mathcal L}, \End({\mathcal M}))$ of cochains of ${\mathcal L}$ with coefficients in $\End({\mathcal M})$ (viewed as a trivial $L_\infty$ module). The DGLA structure on the Chevalley-Eilenberg complex is given by the commutator on $\End({\mathcal M})$ combined with the wedge product. Now, assume that we have an algebra $({\mathcal L,} {\mathcal M})$ over an two-colored operad ${\mathcal P}$ to which the $L_\infty$ operad maps. Then $L_\infty$ ${\mathcal L}$-module structures on ${\mathcal M}$ that are given by universal formulas in terms of operations from ${\mathcal P}$ are the same as Maurer-Cartan elements of the complex $C_{\mathcal P}({\mathcal L}, \End({\mathcal M}))$ of Lie algebra cochains given by universal operations from ${\mathcal P}$.

It is  clear that the $L_\infty$ module structure from Theorem \ref{thm:1} is given by universal formulas in terms of operations from ${\rm{Calc}}_{\rm{alg}}.$ Modulo $u$, these formulas involve only the precalculus analog of ${\rm{Calc}}_{\rm{alg}}.$ Consider three complexes
\begin{equation}\label{eq:univ cochains for calc}
C_{\cal P}^\bullet(\g[\epsilon, u], {\rm{End}}_{K[[u]]}(\Omega^{-\bullet}[[u]]))
\end{equation}
where ${\cal P}$ stands for ${\rm{Calc}}_{\rm{alg}}$, ${\rm{Calc}}_{\infty},$ or ${\rm{Calc}}.$ As explained above, these are complexes of cochains of $\g[\epsilon, u]$ with coefficients in ${\rm{End}}_{K[[u]]}(\Omega^{-\bullet}[[u]])$ that are given by universal operations from ${\rm{Calc}}_{\rm{alg}}$, resp. ${\rm{Calc}}_{\infty},$ resp. ${\rm{Calc}}.$ Here $({\cal{V}}^\bullet , \Omega^\bullet)$ is any algebra over one of the three two-colored operads; $\g$ is ${\cal{V}}^{\bullet+1}$ viewed as an $L_\infty$ algebra via the map of ${\rm{Lie}}_\infty^+$ to one of these operads. In particular, a cochain in (\ref{eq:univ cochains for calc}) produces a cochain in $C^\bullet(\g[\epsilon, u], {\rm{End}}_{K[[u]]}(\Omega^{-\bullet}[[u]]))$ for any ${\cal P}$-algebra.

Recall that a two-colored operad is, in particular, a collection of complexes ${\mathcal{O}}(n)$ and ${\mathcal{M}}(n)$; the first stands for operations ${{\cal{V}}^\bullet}^{\otimes n}\to {\cal{V}}^\bullet;$ the second for operations ${{\cal{V}}^\bullet}^{\otimes n}\otimes \Omega^{\bullet} \to {\Omega}^\bullet.$

First, observe that the three complexes are all quasi-isomorphic. Indeed, they are given by direct sums or products of copies of subspaces of invariants of ${\mathcal{M}}(n)$ with respect to some subgroups of the symmetric group $S_n$, with an extra (Chevalley-Eilenberg) differential. Therefore a quasi-isomorphism of operads leads to a quasi-isomorphism of complexes.

The $L_\infty$ module structures that we are looking for are Maurer-Cartan elements of the DGLAs of cochains of the type as above. The above quasi-isomorphisms preserve the Lie algebra structure and therefore induce isomorphisms on the sets of Maurer-Cartan elements up to equivalence. We want to prove that any two Maurer-Cartan cochains as above, defined by universal operations from $\Ca _{\rm{alg}},$ are equivalent. We see that we can replace $\Ca _{\rm{alg}}$ by $\Ca .$ Therefore, it suffices to prove the following. Let an $L_\infty$ module structure be given by universal formulas purely in terms of the calculus operations $[a,b],\,ab,\,L_a,\,\iota_a, \,d;$ moreover, modulo $u$, it is given by the first four, and it is the original action b) at the level of homology. We claim that any such $L_\infty$ structure is $L_\infty$ equivalent to the original one.

More precisely, we have to prove the following. Let $({\cal{V}}^{\bullet}, \Omega^{\bullet})$ be any DG calculus; $\g=({\cal{V}}^{\bullet+1}, \delta, \{,\});$ consider a Maurer-Cartan element of the DGLA $C^\bullet (\g[\epsilon, u], {\rm{End}}\Omega^{-\bullet}[[u]])$, the cochain complex of $(\g[\epsilon, u], \delta + u\frac{\partial}{\partial \epsilon})$ with coefficients in ${\rm{End}}\Omega^{-\bullet}[[u]]$ on which $\g[\epsilon, u]$ acts trivially. The Maurer-Cartan element, by definition, satisfies
$$(\delta+u\frac{\partial}{\partial \epsilon}+\partial _{\rm{Lie}}){\lambda}+\frac{1}{2}[\lambda,\lambda]=0$$
where $\partial _{\rm{Lie}}$ is the Chevalley-Eilenberg differential. The element $\lambda$ is a cochain defined by universal operations given by formulas involving the five calculus operations; modulo $u$, it involves the four precalculus operations only. For an $n$-linear map $\varphi,$ put
$${\rm{Avg}}\varphi(a_1,\ldots,a_n)=\frac{1}{n!}\sum_{\sigma\in S_n}\pm \varphi(a_{\sigma 1}, \ldots, a_{\sigma _n})$$
The only cochains of suitable degree are of the form
\begin{equation}
\label{eq:coc}
\lambda=\sum_{n\geq 1}\alpha_n \Phi_n+u\sum_{n\geq 1}\beta_n \Psi_n+\sum_{n\geq 0;\,k}\gamma_n^k \Theta_n^k
\end{equation}
where
$$\Phi(a_1\epsilon,\ldots,a_n\epsilon)={\rm{Avg}}L_{a_1}\ldots L_{a_{n-1}}\iota_{a_n};$$
$$\Psi(a_1\epsilon,\ldots,a_n\epsilon)={\rm{Avg}}L_{a_1}\ldots L_{a_{n}}d;$$
$$\Theta_n^k(a_1\epsilon,\ldots,a_n\epsilon, c)={\rm{Avg}}L_{a_1}\ldots L_{a_k}\iota_c  L_{a_{k+1}}\ldots L_{a_n}.$$
First, note that $\gamma_0^0=1$ and $\gamma _n^k=0$ for all $n\geq 1.$ Indeed, the component of $\partial _{\rm{Lie}}$ with values in cochains $\varphi(a_1\epsilon,\ldots,a_n\epsilon, c_1,c_2)$ is zero on all $\Phi_n$ and $\Psi_n;$ for $\Theta_n^k$, it is equal to $\pm {\rm{Avg}}L_{a_1}\ldots L_{a_k}\iota_{\{c_1,c_2\}}  L_{a_{k+1}}\ldots L_{a_n}.$ But this component must be zero for a Maurer-Cartan element.

Note that any cochain $\lambda=\sum_{n\geq 1}\alpha_n \Phi_n$ defines an $L_\infty$ action of $(\g[\epsilon],\delta);$ these actions are non-equivalent. For
$$\lambda=\sum_{n\geq 1}\alpha_n \Phi_n+u\sum_{n\geq 1}\beta_n \Psi_n $$
the terms $\delta \lambda,$ $u\frac{\partial}{\partial\epsilon}\lambda,$  and $\partial _{\rm{Lie}}\lambda$ are all zero. Thus we have
$$ud\lambda+\frac{1}{2}[\lambda,\lambda]=0.$$One has $\beta_1=0$ and $\alpha_1=1$ because the action on the cohomology is the original one. Now, the gauge transformation
$${\rm{exp}}(\sum \kappa_n X_n),$$
$$X_n(a_1\epsilon,\ldots,a_n\epsilon)={\rm{Avg}}L_{a_1}\ldots L_{a_{n-1}}\iota_{a_n}d,$$
kills the $\Psi_n$ terms, $n\geq 2.$ Finally, the Maurer-Cartan equation for the cochain $\sum \alpha_n\Phi_n$ shows that $\alpha_n=0,$ $n\geq 2.$
\subsection{Comparison to Section \ref{s:oofincc}} We conclude this section by citing the result from \cite{TT} that explains the title of Section \ref{s:oofincc}.

 As above, let $A$ be an associative unital algebra over a unital ring $K$ of characteristic zero. Because of Theorem \ref{thm:dtt2}, the pair of complexes $(C^\bullet(A), C_{\bullet}(A))$ is quasi-isomorphic to a pair of complexes $({\mathcal{V}}^\bullet(A), \Omega^\bullet(A))$ that has a structure of a DG calculus.
\begin{teo}\label{thm:tt}
There is an $A_\infty$ quasi-isomorphism
$$(\Y({\mathcal{V}}^\bullet(A))[[u]], \delta+ud)\to {\rm{CC}}^-_{-\bullet}(C^\bullet(A))$$
where the left hand side is an $A_\infty$ algebra as in Theorem \ref{thm: a-infty mod}.
There is also an $A_\infty$ quasi-isomorphism of $A_\infty$ modules
$$(\Omega^{-\bullet}(A)[[u]], b+ud)\to {\rm{CC}}^-_{-\bullet}(A)$$
compatible with the $A_\infty$ map above.
\end{teo}
\section{The Gauss-Manin connection}\label{s:tgmc}
Let ${\mathcal A}$ be a sheaf of ${\mathcal O}_S$-algebras where $S$ is a manifold (real, complex, or algebraic). We assume that ${\mathcal A}$ carries a connection $\nabla$ (not necessarily compatible with the product). Let ${\rm{CC}}^{\rm{per}}_\bullet({\mathcal{A}})$ be the sheaf of periodic cyclic complexes of ${\mathcal A}$ over ${\mathcal O}_S.$ (If ${\mathcal A}$ is the sheaf of local sections of a bundle of algebras then ${\rm{CC}}^{\rm{per}}_\bullet({\mathcal{A}})$ is the sheaf of local sections of the bundle of complexes $s\mapsto {\rm{CC}}^{\rm{per}}_\bullet({\mathcal{A}}_s)$). We conclude the paper by constructing a {\em flat superconnection} on ${\rm{CC}}^{\rm{per}}_\bullet({\mathcal{A}})$, i.e. an operator
$$\nabla_{\rm{GM}}:\Omega^\bullet _S \otimes_{{\mathcal O}_S} {\rm{CC}}^{\rm{per}}_\bullet({\mathcal{A}}) \to \Omega^\bullet _S \otimes_{{\mathcal O}_S} {\rm{CC}}^{\rm{per}}_\bullet({\mathcal{A}})$$
of degree one such that $\nabla_{\rm{GM}}^2=0$ and $\nabla_{\rm{GM}}(fa)=f\nabla_{\rm{GM}}(a)+df\cdot a$ for a function $f$ and a local section $a$.

Let $C^\bullet ({\mathcal A})$ be the sheaf of Hochschild cochain complexes of ${\mathcal A}$ over ${\mathcal O}_S.$ The product on ${\mathcal A}$ defines a two-cochain $m$; then $\nabla m$ is a section of $\Omega ^1(S, C^2 ({\mathcal A})).$ Note also that $\nabla ^2=R\in \Omega ^2(S, {\rm{End}}({\mathcal A}))=\Omega ^2(S,C^0({\mathcal A})).$ Put
$$\alpha = \nabla m+R;$$
one has
$$(\delta+\nabla)^2=\alpha;\;\;\; (\delta+\nabla)(\alpha)=0$$
(recall that $\delta=[m,?]$ is the Hochschild differential).

Now consider an $L_\infty$ module structure as in the Lemma \ref{lemma:relative cochain}. Following a general procedure outlined in \cite{Bar}, Remarks 3.3 and 6.7, put
$$\nabla_{\rm{GM}}=uB+\nabla+\sum_{k,n\geq 0;\;k+n\geq 1} \frac{u^{-n}}{k!n!}\phi_{k+n}(m,\ldots, m, \epsilon \alpha,\ldots,\epsilon \alpha)$$
Here $\phi_{k+n}:S^{k+n}(\g_A[\epsilon, u][1])\to {\rm{End}}({\rm{CC}}^{\rm{per}}_\bullet({\mathcal{A}}))$ are the components of the $L_\infty$ module structure. Note that the summand in $A$ corresponding to $k=1$ and $n=0$ is the Hochschild differential $b$.
\begin{pred} \label{prop:ngm} $\nabla_{\rm{GM}}$ is a flat superconnection.
\end{pred}
{\bf Proof.} Note from the above that the element $\lambda=\nabla+m+u^{-1}\epsilon\alpha$ satisfies $[\lambda,\lambda]=0$ in the DG Lie algebra ${\mathcal L}^\bullet [\epsilon, u]$ where ${\mathcal L}$ is the DGLA of Hochschild cochains of the graded algebra 
$$\Omega^\bullet _S({\mathcal A})=\Omega^\bullet_S\otimes_{{\mathcal {O}}_S}{\mathcal A}.$$ 
Therefore the operator
 $$uB+\sum _{n=1}^\infty \frac{1}{n!}\phi_n(\lambda,
\ldots, \lambda)$$
on the periodic cyclic complex of $\Omega^\bullet _S({\mathcal A})$ has square zero. Because $\phi_n(\nabla, \ldots)=0$ for $n>1,$ this operator descends to the quotient ${\rm{CC}}^{\rm{per}}_\bullet({\mathcal{A}})$; the resulting operator is exactly $\nabla _{\rm{GM}}$.

Of course the flatness of $\nabla _{\rm{GM}}$ can be obtained directly from the $L_\infty$ identities. In fact, if we denote the infinite sum in the above formula by $A$ and write $\nabla_{\rm{GM}}=uB+\nabla+A$, we have
$$(\nabla_{\rm{GM}})^2 = [uB, A]+[\nabla, A]+\nabla ^2+\frac{1}{2}[A,A]=$$
$$-\sum_{k,n\geq 0;\;k+n\geq 1} \frac{u^{-n}}{k!(n-1)!}\phi_{k+n}(m,\ldots, m, \epsilon \delta\alpha, \epsilon\alpha, \ldots,\epsilon \alpha)$$
$$-\sum_{k,n\geq 0;\;k+n\geq 1} \frac{u^{-n+1}}{k!(n-1)!}\phi_{k+n}(m,\ldots, m, \alpha, \epsilon\alpha, \ldots,\epsilon \alpha)-\frac{1}{2}[A,A]+$$
$$\sum_{k,n\geq 0;\;k+n\geq 1} \frac{u^{-n}}{(k-1)!n!}\phi_{k+n}(\nabla m, m,\ldots, m, \epsilon \alpha,\ldots,\epsilon \alpha)+$$
$$-\sum_{k,n\geq 0;\;k+n\geq 1} \frac{u^{-n}}{k!(n-1)!}\phi_{k+n}(m,\ldots, m, \epsilon \nabla\alpha, \epsilon\alpha, \ldots,\epsilon \alpha)+R+\frac{1}{2}[A,A]=0$$
Here we used the facts that $[m,m]=0,$ $[\epsilon\alpha,\epsilon\alpha]=0,$ $(\delta+\nabla)\alpha=0,$ $[m,\epsilon\alpha]=-\epsilon\delta\alpha,$ and $\phi_n(R, \ldots)=0$ for $n>1.$

 ~\\

\noindent\textsc{Department of Mathematics,
University of California at Riverside, \\
900 Big Springs Drive,\\
Riverside, CA 92521, USA \\
\emph{E-mail address:} {\bf vald@math.ucr.edu}}

~\\

\noindent\textsc{Mathematics Department,
Northwestern University, \\
2033 Sheridan Rd.,\\
Evanston, IL 60208, USA \\
\emph{E-mail addresses:} {\bf tamarkin@math.northwestern.edu},
{\bf tsygan@math.northwestern.edu}}


\begin{thebibliography}{99}
\bibitem{Bar} S.~Barannikov, {\emph Quantum periods. I. Semi-infinite variations of Hodge structures}, Internat. Math. Res. Notices  2001,  no. 23, 1243--1264.
\bibitem{Baran} V. Baranovsky, {\em A universal enveloping of an $L_{\infty}$ algebra}, arxiv:0706.1396
\bibitem{bgnt1}
P.~Bressler, A.~Gorokhovsky, R.~Nest, and B.~Tsygan, {\em Algebraic index theorem for symplectic deformations of gerbes}, to appear.
\bibitem{BNT} P.~Bressler, R.~Nest, R., B.~Tsygan, \emph{Riemann-Roch theorems via deformation quantization, I, II},  Adv. Math.  \textbf{167}  (2002),  no. 1, 1--25, 26--73.
\bibitem{CR} D. Calaque and C. A. Rossi,
{\em Shoikhet's Conjecture and Duflo Isomorphism on (Co)Invariants},
SIGMA {\bf 4} (2008) 060; arXiv:0805.2409.
 \bibitem{Ca1} A. Caldararu,  {\em The Mukai pairing, I: the Hochschild structure}, arXiv:math/0308079.
\bibitem{Ca2} A. Caldararu, {\em The Mukai pairing, II: the Hochschild-Kostant-Rosenberg isomorphism}, arXiv:math/0308080.
\bibitem{Ca3} A. Caldararu, S. Willerton, {\em The Mukai pairing, I: a categorical approach}, arXiv:0707.2052.
\bibitem{CE} H.~ Cartan, S.~ Eilenberg, {\em Homological Algebra}, Princeton, 1956.
.
\bibitem{Co} K. Costello, {\em Topological conformal field theories and Calabi-Yau categories}, arXiv:math/0605647.
\bibitem{D} V.A. Dolgushev,
{\em A Proof of Tsygan's formality conjecture
for an arbitrary smooth manifold}, PhD thesis,
MIT; math.QA/0504420.
\bibitem{DGT} Yu. Daletski, I. Gelfand, and
B. Tsygan, {\em On a variant of noncommutative
geometry}, Soviet Math. Dokl. {\bf 40}, 2
(1990) 422--426.
\bibitem{DG} I.Ya. Dorfman, I. M. Gelfand, {\em Hamiltonian operators and algebraic structures related to them,} Funct. Anal. Appl. {\bf 13} (1979), pp. 13–30.
\bibitem{DTT1} V. Dolgushev, D. Tamarkin and
B. Tsygan, {\em The homotopy Gerstenhaber algebra of Hochschild
cochains of a regular algebra is formal}, J. Noncommut. Geom.
{\bf 1}, 1 (2007) 1--25;
arXiv:math/0605141.
\bibitem{DTT2}  V. Dolgushev, D. Tamarkin, and
B. Tsygan, {\em Formality of the homotopy calculus algebra of
Hochschild (co)chains}, arXiv:0807.5117.
\bibitem{DTT3}  V. Dolgushev, D. Tamarkin, and
B. Tsygan, {\em Formality theorems for Hochschild complexes and their applications}, Submitted to: Proceedings of the Conference Poisson 2008, Lausanne.
\bibitem{Ge} M.~Gerstenhaber, {\em The cohomology structure of an associative ring}, Ann. Math. (2), {\bf 78} (1963), 267--288.
\bibitem{G} E. Getzler, {\em Cartan homotopy formulas and the Gauss-Manin connection in cyclic homology}, Quantum deformations of algebras and their representations (Ramat-Gan, 1991/1992; Rehovot, 1991/1992), 65--78, Israel Math. Conf. Proc., 7, Bar-Ilan Univ., Ramat Gan, 1993.
\bibitem{GJ} E.~Getzler and J.~Jones, {\it $A_{\infty}$ algebras and the cyclic bar complex}, Illinois J. of Math. {\bf 34} (1990), 256-283.
\bibitem{GJ1} E.~Getzler,~J. Jones, {\em Operads, homotopy algebra
    and iterated integrals for double loop spaces}, preprint hep-th9403055.
\bibitem{Gin} V. Ginzburg, T. Schedler, {\em Free products, cyclic homology, and the Gauss-Manin connection,} arXiv:0803.3655.
\bibitem{Hin} V. Hinich, {\em Tamarkin's proof of
Kontsevich formality theorem}, Forum Math. {\bf 15},
4 (2003) 591--614; math.QA/0003052.
\bibitem{HJ} C.E.~Hood and J.D.S.~Jones, {\it Some algebraic properties of cyclic homology groups}, K - Theory, {\bf 1} (1987), 361-384.
\bibitem{Ka} T. V. Kadeishvili, {\em On the homology theory of fibered spaces}, Uspekhi Mat. Nauk, {\bf 35} (3), 1980, 183-188.
    \bibitem{KKP} L.Katzarkov, M.Kontsevich,
T.Pantev, {\em Hodge theoretic aspects of mirror symmetry}, arXiv:0806.0107.
\bibitem{Kel} B. Keller, {\em $A_\infty$ algebras, modules and functor categories}, Trends in representation theory of algebras and related topics, 67--93, Contemporary Mathematics, \textbf{406}, AMS, 2006.
\bibitem{K1} M. Kontsevich, {\em Operads and motives in deformation quantization}, Lett. Math. Phys., {\bf 48} (1999) 35--72.
\bibitem{KS} M. Kontsevich, Y. Soibelman, {\em Notes on A-infinity algebras, A-infinity categories and non-commutative geometry. I}, math.RA/0606241.
\bibitem{KS1} M.~Kontsevich, Y.~Soibelman, {\em Deformations of algebras over operads and the Deligne conjecture}, Conf\'{e}rence Moshe Flato 1999, Vol. 1 (Dijon), 255-307, Math. Phys. Studies {\bf 21}, Kl\"{u}wer Academic Publishers, Dordrecht, 2000.
\bibitem{LS} T. Lada and J. Stasheff, {\em Introduction to SH Lie
algebras for physicists}, Intern. J. Theor Phys. {\bf 32}, 7 (1993)
1087-1103.
\bibitem{Loday} J.- L. Loday, {\em Cyclic Homology},
Gr\"{u}ndlehren der mathematischen Wissenschaften, 301.
Springer-Verlag, Berlin, 1992.
\bibitem{Ma} N. Markarian, {\em The Atiyah class, Hochschild cohomology and the Riemann-Roch theorem}, arXiv:math/0610553.
\bibitem{Me} S. Merkulov, {\em An $L_\infty$ algebra of an unobstructed deformation functor}, Int. Math. Res. Notices (2000), no.3, 147--164.
\bibitem{Me1} S. Merkulov, {\em Exotic automorphisms of the Schouten algebra of polyvector fields},  arXiv:0809.2385.
\bibitem{MV} S. Merkulov, B. Vallette,
{\em Deformation theory of representations of prop(erad)s},  arXiv:0707.0889.
\bibitem{NT} R. Nest, B. Tsygan, {\em On the cohomology ring of an algebra}, in: Advances in geometry, Progress in Mathematics, {\bf 172}, Birkh\"{a}user, 337--370.
\bibitem{Ra} A. Ramadoss, {\em The Mukai pairing and integral transforms in Hochschild homology}, arXiv:0805.1760.
\bibitem{Sh} B. Shoikhet,
{\em A proof of the Tsygan formality conjecture for chains},
Adv. Math., {\bf 179}, 1 (2003) 7--37;
math.QA/0010321.
\bibitem{St} J.~ Stasheff, {\em Homotopy associativity of H-spaces, I and
    II}, Trans. AMS {\bf 108} (1963), 275-312.
\bibitem{Ta1} D. Tamarkin,
{\em Another proof of M. Kontsevich formality theorem},
math.QA/9803025.
\bibitem{Tam} D. Tamarkin,
{\em Formality of chain operad of little discs},
Lett. Math. Phys.  {\bf 66}, 1-2  (2003) 65--72;
math.QA/9809164.
\bibitem{Ta} D. Tamarkin, {\em What do DG categories
form?}  Compos. Math. {\bf 143}, 5 (2007) 1335--1358;
math.CT/0606553.
\bibitem{TT} D. Tamarkin and B. Tsygan, {\em The ring of differential operators on forms in noncommutative calculus}, Graph patterns in mathematics and theoretical physics, 105--131, Proc. Symp. Pure Math. {\bf 73}, AMS, Providence, RI, 2005.
\bibitem{TT1} D. Tamarkin and B. Tsygan,
{\em Cyclic formality and index theorems},
Talk given at the Mosh\'e Flato Conference (2000),
Lett. Math. Phys. {\bf 56}, 2  (2001) 85--97.
\bibitem{ToV} B. To\"{e}n, G. Vezzosi,
{\em A note on Chern character, loop spaces and derived algebraic geometry}, Proc. of the 2007 Abel Symposium, arXiv:0804.1274.
\bibitem{T} B. Tsygan, {\em On the Gauss-Manin connection in cyclic homology}, Methods Funkt. Anal. Topology, {\bf 13}, 1, 83--94, 2007.
\end{thebibliography}
\end{document}